# ON THE DISTRIBUTION OF THE MAXIMUM OF A GAUSSIAN FIELD WITH $D$ PARAMETERS[1]

By Jean-Marc Azaïs and Mario Wschebor

*Université Paul Sabatier and Universidad de la República*

Let $I$ be a compact $d$-dimensional manifold, let $X : I \to \mathcal{R}$ be a Gaussian process with regular paths and let $F_I(u)$, $u \in \mathcal{R}$, be the probability distribution function of $\sup_{t \in I} X(t)$.

We prove that under certain regularity and nondegeneracy conditions, $F_I$ is a $C^1$-function and satisfies a certain implicit equation that permits to give bounds for its values and to compute its asymptotic behavior as $u \to +\infty$. This is a partial extension of previous results by the authors in the case $d = 1$.

Our methods use strongly the so-called Rice formulae for the moments of the number of roots of an equation of the form $Z(t) = x$, where $Z : I \to \mathcal{R}^d$ is a random field and $x$ is a fixed point in $\mathcal{R}^d$. We also give proofs for this kind of formulae, which have their own interest beyond the present application.

**1. Introduction and notation.** Let $I$ be a $d$-dimensional compact manifold and let $X : I \to \mathcal{R}$ be a Gaussian process with regular paths defined on some probability space $(\Omega, \mathcal{A}, \mathrm{P})$. Define $M_I = \sup_{t \in I} X(t)$ and $F_I(u) = \mathrm{P}\{M_I \le u\}$, $u \in \mathcal{R}$, the probability distribution function of the random variable $M_I$. Our aim is to study the regularity of the function $F_I$ when $d > 1$.

There exist a certain number of general results on this subject, starting from the papers by Ylvisaker (1968) and Tsirelson (1975) [see also Weber (1985), Lifshits (1995), Diebolt and Posse (1996) and references therein]. The main purpose of this paper is to extend to $d > 1$ some of the results about the regularity of the function $u \rightsquigarrow F_I(u)$ in Azaïs and Wschebor (2001), which concern the case $d = 1$.

Our main tool here is the Rice formula for the moments of the number of roots $N_u^Z(I)$ of the equation $Z(t) = u$ on the set $I$, where $\{Z(t) : t \in I\}$ is an

Received January 2003; revised November 2003.

[1]Supported by ECOS program U97E02.

$\mathcal{R}^d$-valued Gaussian field, $I$ is a subset of $\mathcal{R}^d$ and $u$ is a given point in $\mathcal{R}^d$. For $d > 1$, even though it has been used in various contexts, as far as the authors know, a full proof of the Rice formula for the moments of $N_u^Z(I)$ seems to have only been published by Adler (1981) for the first moment of the number of critical points of a real-valued stationary Gaussian process with a $d$-dimensional parameter, and extended by Azaïs and Delmas (2002) to the case of processes with constant variance. Cabaña (1985) contains related formulae for random fields; see also the Ph.D. thesis of Konakov cited by Piterbarg (1996b). In the next section we give a more general result which has an interest that goes beyond the application of the present paper. At the same time the proof appears to be simpler than previous ones. We have also included the proof of the formula for higher moments, which in fact follows easily from the first moment. Both extend with no difficulties to certain classes of non-Gaussian processes.

It should be pointed out that the validity of the Rice formula for Lebesgue-almost every $u \in \mathcal{R}^d$ is easy to prove [Brillinger (1972)] but this is insufficient for a certain number of standard applications. For example, assume $X : I \rightsquigarrow \mathcal{R}$ is a real-valued random process and one is willing to compute the moments of the number of critical points of $X$. Then, we must take for $Z$ the random field $Z(t) = X'(t)$ and the formula one needs is for the precise value $u = 0$ so that a formula for almost every $u$ does not solve the problem.

We have added the Rice formula for processes defined on smooth manifolds. Even though the Rice formula is local, this is convenient for various applications. We will need a formula of this sort to state and prove the implicit formulae for the derivatives of the distribution of the maximum (see Section 3).

The results on the differentiation of $F_I$ are partial extensions of Azaïs and Wschebor (2001). Here, we have only considered the first derivative $F_I'(u)$. In fact, one can push our procedure one step more and prove the existence of $F_I''(u)$ as well as some implicit formula for it. But we have not included this in the present paper since formulae become very complicated and it is unclear at present whether the actual computations can be performed, even in simple examples. The technical reason for this is that, following the present method, to compute $F_I''(u)$, one needs to differentiate expressions that contain the "helix process" that we introduce in Section 4, containing singularities with unpleasant behavior [see Azaïs and Wschebor (2002)].

For Gaussian fields defined on a $d$-dimensional regular manifold ($d > 1$) and possessing regular paths we obtain some improvements with respect to classical and general results due to Tsirelson (1975) for Gaussian sequences. An example is Corollary 5.1, which provides an asymptotic formula for $F_I'(u)$ as $u \to +\infty$ which is explicit in terms of the covariance of the process and can be compared with Theorem 4 in Tsirelson (1975) where an implicit expression depending on the function $F$ itself is given.

We use the following notation:

If $Z$ is a smooth function $U \rightsquigarrow \mathcal{R}^{d'}$, $U$ a subset of $\mathcal{R}^d$, its successive derivatives are denoted $Z'$, $Z''$, ..., $Z^{(k)}$ and considered, respectively, as linear, bilinear, ..., $k$-linear forms on $\mathcal{R}^d$. For example, $X^{(3)}(t)[v_1, v_2, v_3] = \sum_{i,j,k=1}^d \frac{\partial^{3X(t)}}{\partial t^i \partial t^j \partial t^k} v_1^i v_2^j v_3^k$. The same notation is used for a derivative on a $\mathcal{C}^\infty$ manifold.

$\dot{I}, \partial I$ and $\bar{I}$ are, respectively, the interior, the boundary and the closure of the set $I$. If $\xi$ is a random vector with values in $\mathcal{R}^d$, whenever they exist, we denote by $p_\xi(x)$ the value of the density of $\xi$ at the point $x$, by $\mathrm{E}(\xi)$ its expectation and by $\mathrm{Var}(\xi)$ its variance–covariance matrix. $\lambda$ is Lebesgue measure. If $u, v$ are points in $\mathcal{R}^d$, $\langle u, v \rangle$ denotes their usual scalar product and $\|u\|$ the Euclidean norm of $u$. For $M$ a $d \times d$ real matrix, we denote $\|M\| = \sup_{\|x\|=1} \|Mx\|$.

Also for symmetric $M$, $M \succ 0$ (resp. $M \prec 0$) denotes that $M$ is positive definite (resp. negative definite). $A^c$ denotes the complement of the set $A$. For real $x$, $x^+ = \sup(x, 0)$, $x^- = \sup(-x, 0)$.

**2. Rice formulae.** Our main results in this section are the following:

THEOREM 2.1. *Let* $Z : I \rightsquigarrow \mathcal{R}^d$, $I$ *a compact subset of* $\mathcal{R}^d$, *be a random field and* $u \in \mathcal{R}^d$. *Assume that:*

(A0) *$Z$ is Gaussian.*
(A1) *$t \rightsquigarrow Z(t)$ is a.s. of class $\mathcal{C}^1$.*
(A2) *For each $t \in I$, $Z(t)$ has a nondegenerate distribution [i.e., $\mathrm{Var}(Z(t)) \succ 0$].*
(A3) *$\mathrm{P}\{\exists t \in \dot{I}, Z(t) = u, \det(Z'(t)) = 0\} = 0$.*
(A4) *$\lambda(\partial I) = 0$.*

*Then*

$$(1) \qquad \mathrm{E}(N_u^Z(I)) = \int_I \mathrm{E}(|\det(Z'(t))| / Z(t) = u) p_{Z(t)}(u) \, dt,$$

*and both members are finite.*

THEOREM 2.2. *Let $k$, $k \geq 2$, be an integer. Assume the same hypotheses as in Theorem 2.1 except for* (A2), *which is replaced by:*

(A'2) *for $t_1, \ldots, t_k \in I$ pairwise different values of the parameter, the distribution of $(Z(t_1), \ldots, Z(t_k))$ does not degenerate in $(\mathcal{R}^d)^k$. Then*

$$\mathrm{E}[(N_u^Z(I))(N_u^Z(I) - 1) \cdots (N_u^Z(I) - k + 1)]$$

$$(2) \qquad = \int_{I^k} \mathrm{E}\left( \prod_{j=1}^k |\det(Z'(t_j))| / Z(t_1) = \cdots = Z(t_k) = u \right)$$

$$\times \, p_{Z(t_1), \ldots, Z(t_k)}(u, \ldots, u) \, dt_1 \cdots dt_k,$$

*where both members may be infinite.*

REMARK. Note that Theorem 2.1 (resp. Theorem 2.2) remains valid if one replaces $I$ by $\dot{I}$ in (1) or (2) and if hypotheses (A0)–(A2) [resp. (A'2)] and (A3) are verified. This follows immediately from the above statements. A standard extension argument shows that (1) holds true if one replaces $I$ by any Borel subset of $\dot{I}$.

Sufficient conditions for hypothesis (A3) to hold are given by the next proposition. Under condition (a) the result is proved in Lemma 5 of Cucker and Wschebor (2003). Under condition (b) the proof is straightforward.

PROPOSITION 2.1. *Let* $Z : I \rightsquigarrow \mathcal{R}^d$, *I a compact subset of* $\mathcal{R}^d$, *be a random field with paths of class* $\mathcal{C}^1$ *and* $u \in \mathcal{R}^d$. *Assume that:*

(i) $p_{Z(t)}(x) \leq C$ *for all* $t \in I$ *and* $x$ *in some neighborhood of* $u$.
(ii) *At least one of the two following hypotheses is satisfied:*

    (a) *a.s.* $t \rightsquigarrow Z(t)$ *is of class* $\mathcal{C}^2$,
    (b) $\alpha(\delta) = \sup_{t \in I, x \in V(u)} \mathrm{P}\{|\det(Z'(t))| < \delta / Z(t) = x\} \to 0$ *as* $\delta \to 0$,
*where* $V(u)$ *is some neighborhood of* $u$.

*Then* (A3) *holds true.*

The following lemma is easy to prove.

LEMMA 2.1. *With the notation of Theorem 2.1, suppose that* (A1) *and* (A4) *hold true and that* $p_{Z(t)}(x) \leq C$ *for all* $t \in I$ *and* $x$ *in some neighborhood of* $u$. *Then* $\mathrm{P}\{N_u^Z(\partial I) \neq 0\} = 0$.

LEMMA 2.2. *Let* $Z : I \to \mathcal{R}^d$, *I a compact subset of* $\mathcal{R}^d$, *be a* $\mathcal{C}^1$ *function and* $u$ *a point in* $\mathcal{R}^d$. *Assume that:*

(a) $\inf_{t \in Z^{-1}(\{u\})}(\lambda_{\min}(Z'(t))) \geq \Delta > 0$,
(b) $\omega_{Z'}(\eta) < \Delta/d$,

*where* $\omega_{Z'}$ *is the continuity modulus of* $Z'$, *defined as the maximum of the continuity moduli of its entries,* $\lambda_{\min}(M)$ *is the square root of the smallest eigenvalue of* $M^T M$ *and* $\eta$ *is a positive number.*

*Then, if* $t_1, t_2$ *are two distinct roots of the equation* $Z(t) = u$ *such that the segment* $[t_1, t_2]$ *is contained in* $I$, *the Euclidean distance between* $t_1$ *and* $t_2$ *is greater than* $\eta$.

PROOF. Set $\tilde{\eta} = \|t_1 - t_2\|$, $v = \frac{t_1 - t_2}{\|t_1 - t_2\|}$. Using the mean value theorem, for $i = 1, \ldots, d$, there exists $\xi_i \in [t_1, t_2]$ such that $(Z'(\xi_i)v)_i = 0$. Thus

$$|(Z'(t_1)v)_i| = |(Z'(t_1)v)_i - (Z'(\xi_i)v)_i|$$

$$\leq \sum_{k=1}^{d} |Z'(t_1)_{ik} - Z'(\xi_i)_{ik}| |v_k| \leq \omega_{Z'}(\tilde{\eta}) \sum_{k=1}^{d} |v_k| \leq \omega_{Z'}(\tilde{\eta})\sqrt{d}.$$

In conclusion, $\Delta \leq \lambda_{\min}(Z'(t_1)) \leq \|Z'(t_1)v\| \leq \omega_{Z'}(\tilde{\eta})d$, which implies $\tilde{\eta} > \eta$. □

PROOF OF THEOREM 2.1. Consider a continuous nondecreasing function $F$ such that $F(x) = 0$ for $x \leq 1/2$, $F(x) = 1$ for $x \geq 1$. Let $\Delta$ and $\eta$ be positive real numbers. Define the random function

$$\alpha_{\Delta, \eta}(u) = F\left(\frac{1}{2\Delta} \inf_{s \in I} [\lambda_{\min}(Z'(s)) + \|Z(s) - u\|]\right) \times \left(1 - F\left(\frac{d}{\Delta}\omega_{Z'}(\eta)\right)\right),$$
(3)

and the set $I_{-\eta} = \{t \in I : \|t - s\| \geq \eta, \forall s \notin I\}$. If $\alpha_{\Delta, \eta}(u) > 0$ and $N_u^Z(I_{-\eta})$ does not vanish, conditions (a) and (b) in Lemma 2.2 are satisfied. Hence, in each ball with diameter $\frac{\eta}{2}$ centered at a point in $I_{-\eta}$, there is at most one root of the equation $Z(t) = u$, and a compactness argument shows that $N_u^Z(I_{-\eta})$ is bounded by a constant $C(\eta, I)$, depending only on $\eta$ and on the set $I$.

Take now any real-valued nonrandom continuous function $f : \mathcal{R}^d \to \mathcal{R}$ with compact support. Because of the coarea formula [Federer (1969), Theorem 3.2.3], since a.s. $Z$ is Lipschitz and $\alpha_{\Delta, \eta}(u) \cdot f(u)$ is integrable,

$$\int_{\mathcal{R}^d} f(u) N_u^Z(I_{-\eta}) \alpha_{\Delta, \eta}(u) \, du = \int_{I_{-\eta}} |\det(Z'(t))| f(Z(t)) \alpha_{\Delta, \eta}(Z(t)) \, dt.$$

Taking expectations in both sides,

$$\int_{\mathcal{R}^d} f(u) \mathrm{E}(N_u^Z(I_{-\eta}) \alpha_{\Delta, \eta}(u)) \, du$$

$$= \int_{\mathcal{R}^d} f(u) \, du \int_{I_{-\eta}} \mathrm{E}(|\det(Z'(t))| \alpha_{\Delta, \eta}(u)/Z(t) = u) p_{Z(t)}(u) \, dt.$$

It follows that the two functions

(i) $\mathrm{E}(N_u^Z(I_{-\eta}) \alpha_{\Delta, \eta}(u))$,
(ii) $\int_{I_{-\eta}} \mathrm{E}(|\det(Z'(t))| \alpha_{\Delta, \eta}(u)/Z(t) = u) p_{Z(t)}(u) \, dt$,

coincide Lebesgue-almost everywhere as functions of $u$.

Let us prove that both functions are continuous, hence they are equal for every $u \in \mathcal{R}^d$.

Fix $u = u_0$ and let us show that the function in (i) is continuous at $u = u_0$. Consider the random variable inside the expectation sign in (i). Almost surely, there is no point $t$ in $Z^{-1}(\{u_0\})$ such that $\det(Z'(t)) = 0$. By the local inversion theorem, $Z(\cdot)$ is invertible in some neighborhood of each point belonging to $Z^{-1}(\{u_0\})$ and the distance from $Z(t)$ to $u_0$ is bounded below by a positive number for $t \in I_{-\eta}$ outside of the union of these neighborhoods. This implies that, a.s., as a function of $u$, $N_u^Z(I_{-\eta})$ is constant in some (random) neighborhood of $u_0$. On the other hand, it is clear from its definition that the function $u \rightsquigarrow \alpha_{\Delta,\eta}(u)$ is continuous and bounded. We may now apply dominated convergence as $u \to u_0$, since $N_u^Z(I_{-\eta}) \alpha_{\Delta,\eta}(u)$ is bounded by a constant that does not depend on $u$.

For the continuity of (ii), it is enough to prove that, for each $t \in I$ the conditional expectation in the integrand is a continuous function of $u$. Note that the random variable $|\det(Z'(t))| \alpha_{\Delta,\eta}(u)$ is a functional defined on $\{(Z(s), Z'(s)) : s \in I\}$. Perform a Gaussian regression of $(Z(s), Z'(s)) : s \in I$ with respect to the random variable $Z(t)$, that is, write

$$Z(s) = Y^t(s) + \alpha^t(s) Z(t),$$
$$Z'_j(s) = Y^t_j(s) + \beta^t_j(s) Z(t), \qquad j = 1, \ldots, d,$$

where $Z'_j(s)$, $j = 1, \ldots, d$, denote the columns of $Z'(s)$, $Y^t(s)$ and $Y^t_j(s)$ are Gaussian vectors, independent of $Z(t)$ for each $s \in I$, and the regression matrices $\alpha^t(s)$, $\beta^t_j(s)$, $j = 1, \ldots, d$, are continuous functions of $s, t$ [take into account (A2)]. Replacing in the conditional expectation, we are now able to get rid of the conditioning, and using the fact that the moments of the supremum of an a.s. bounded Gaussian process are finite, the continuity in $u$ follows by dominated convergence.

So, now we fix $u \in \mathcal{R}^d$ and make $\eta \downarrow 0$, $\Delta \downarrow 0$ in that order, both in (i) and (ii). For (i) one can use Beppo Levi's theorem. Note that almost surely $N_u^Z(I_{-\eta}) \uparrow N_u^Z(\dot{I}) = N_u^Z(I)$, where the last equality follows from Lemma 2.1. On the other hand, the same Lemma 2.1 plus (A3) imply together that, almost surely,

$$\inf_{s \in I}[\lambda_{\min}(Z'(s)) + \|Z(s) - u\|] > 0$$

so that the first factor in the right-hand side of (3) increases to 1 as $\Delta$ decreases to zero. Hence by Beppo Levi's theorem,

$$\lim_{\Delta \downarrow 0} \lim_{\eta \downarrow 0} \mathrm{E}(N_u^Z(I_{-\eta}) \alpha_{\Delta,\eta}(u)) = \mathrm{E}(N_u^Z(I)).$$

For (ii), one can proceed in a similar way after deconditioning obtaining (1). To finish the proof, remark that standard Gaussian calculations show the finiteness of the right-hand side of (1).   $\square$

PROOF OF THEOREM 2.2. For each $\delta > 0$, define the domain

$$D_{k,\delta}(I) = \{(t_1, \ldots, t_k) \in I^k, \|t_i - t_j\| \geq \delta \text{ if } i \neq j, i,j = 1, \ldots, k\}$$

and the process $\widetilde{Z}$

$$(t_1, \ldots, t_k) \in D_{k,\delta}(I) \rightsquigarrow \widetilde{Z}(t_1, \ldots, t_k) = (Z(t_1), \ldots, Z(t_k)).$$

It is clear that $\widetilde{Z}$ satisfies the hypotheses of Theorem 2.1 for every value $(u, \ldots, u) \in (\mathcal{R}^d)^k$. So,

$$\mathrm{E}[N^{\widetilde{Z}}_{(u,\ldots,u)}(D_{k,\delta}(I))]$$

$$(4) \qquad = \int_{D_{k,\delta}(I)} \mathrm{E}\left(\prod_{j=1}^{k} |\det(Z'(t_j))| / Z(t_1) = \cdots = Z(t_k) = u\right)$$

$$\times \, p_{Z(t_1),\ldots,Z(t_k)}(u, \ldots, u) \, dt_1 \cdots dt_k.$$

To finish, let $\delta \downarrow 0$, note that $(N^Z_u(I))(N^Z_u(I) - 1) \ldots (N^Z_u(I) - k + 1)$ is the monotone limit of $N^{\widetilde{Z}}_{(u,\ldots,u)}(D_{k,\delta}(I))$, and that the diagonal $D_k(I) = \{(t_1, \ldots, t_k) \in I^k, \, t_i = t_j \text{ for some pair } i,j, i \neq j\}$ has zero Lebesgue measure in $(\mathcal{R}^d)^k$. □

REMARK. Even thought we will not use this in the present paper, we point out that it is easy to adapt the proofs of Theorems 2.1 and 2.2 to certain classes of non-Gaussian processes.

For example, the statement of Theorem 2.1 remains valid if one replaces hypotheses (A0) and (A2), respectively, by the following (B0) and (B2):

(B0) $Z(t) = H(Y(t))$ for $t \in I$, where $Y : I \to \mathcal{R}^n$ is a Gaussian process with $\mathcal{C}^1$ paths such that for each $t \in I, Y(t)$ has a nondegenerate distribution and $H : \mathcal{R}^n \to \mathcal{R}^d$ is a $\mathcal{C}^1$ function.

(B2) For each $t \in I, Z(t)$ has a density $p_{Z(t)}$ which is continuous as a function of $(t, u)$.

Note that (B0) and (B2) together imply that $n \geq d$. The only change to be introduced in the proof of the theorem is in the continuity of (ii) where the regression is performed on $Y(t)$ instead of $Z(t)$.

Similarly, the statement of Theorem 2.2 remains valid if we replace (A0) by (B0) and add the requirement that the joint density of $Z(t_1), \ldots, Z(t_k)$ be a continuous function of $t_1, \ldots, t_k, u$ for pairwise different $t_1, \ldots, t_k$.

Now consider a process $X$ from $I$ to $\mathcal{R}$ and define

$$M^X_{u,1}(I) = \sharp\{t \in I, X(\cdot) \text{ has a local maximum at the point } t, X(t) > u\},$$

$$M^X_{u,2}(I) = \sharp\{t \in I, X'(t) = 0, X(t) > u\}.$$

The problem of writing Rice formulae for the factorial moments of these random variables can be considered as a particular case of the previous one and the proofs are the same, mutatis mutandis. For further use, we state as a theorem the Rice formula for the expectation. For breavity we do not state the equivalent of Theorem 2.2, which holds true similarly.

THEOREM 2.3. *Let* $X : I \rightsquigarrow \mathcal{R}$, *I a compact subset of* $\mathcal{R}^d$, *be a random field. Let* $u \in \mathcal{R}$, *define* $M_{u,i}^X(I)$, $i = 1, 2$, *as above. For each* $d \times d$ *real symmetric matrix* $M$, *we put* $\delta^1(M) := |\det(M)| \mathbb{1}_{M \prec 0}$, $\delta^2(M) := |\det(M)|$.

*Assume:*

(A0)  $X$ *is Gaussian,*

(A''1)  *a.s.* $t \rightsquigarrow X(t)$ *is of class* $\mathcal{C}^2$,

(A''2)  *for each* $t \in I$, $X(t), X'(t)$ *has a nondegenerate distribution in* $\mathcal{R}^1 \times \mathcal{R}^d$,

(A''3)  *either a.s.* $t \rightsquigarrow X(t)$ *is of class* $\mathcal{C}^3$ *or* $\alpha(\delta) = \sup_{t \in I, x' \in V(0)} \mathrm{P}(|\det(X''(t))| < \delta / X'(t) = x') \to 0$ *as* $\delta \to 0$, *where* $V(0)$ *denotes some neighborhood of* $0$,

(A4)  $\partial I$ *has zero Lebesgue measure.*

*Then, for* $i = 1, 2$,

$$\mathrm{E}(M_{u,i}^X(I)) = \int_u^\infty dx \int_I \mathrm{E}(\delta^i(X''(t))/X(t) = x, X'(t) = 0) p_{X(t),X'(t)}(x,0) \, dt$$

*and both members are finite.*

2.1. *Processes defined on a smooth manifold.* Let $U$ be a differentiable manifold (by differentiable we mean infinitely differentiable) of dimension $d$. We suppose that $U$ is orientable in the sense that there exists a nonvanishing differentiable $d$-form $\Omega$ on $U$. This is equivalent to assuming that there exists an atlas $((U_i, \phi_i); i \in I)$ such that for any pair of intersecting charts $(U_i, \phi_i)$, $(U_j, \phi_j)$, the Jacobian of the map $\phi_i \circ \phi_j^{-1}$ is positive.

We consider a Gaussian stochastic process with real values and $\mathcal{C}^2$ paths $X = \{X(t) : t \in U\}$ defined on the manifold $U$. In this section we first write Rice formulae for this kind of processes without further hypotheses on $U$. When $U$ is equipped with a Riemannian metric, we give, without details and proof, a nicer form. Other forms exist also when $U$ is naturally embedded in a Euclidean space, but we do not need this in the sequel [see Azaïs and Wschebor (2002)].

We will assume that in every chart $X(t)$ and $DX(t)$ have a nondegenerate joint distribution and that hypothesis (A''3) is verified. For $S$ a Borel subset of $\dot{U}$, the following quantities are well defined and measurable: $M_{u,1}^X(S)$, the number of local maxima and $M_{u,2}^X(S)$, the number of critical points.

PROPOSITION 2.2. *For* $k = 1, 2$ *the quantity which is expressed in every chart* $\phi$ *with coordinates* $s_1, \ldots, s_d$ *as*

$$(5) \quad \int_u^{+\infty} dx \, \mathrm{E}(\delta^k(Y''(s))/Y(s) = x, Y'(s) = 0) p_{Y(s), Y'(s)}(x, 0) \bigwedge_{i=1}^d ds_i,$$

*where* $Y(s)$ *is the process* $X$ *written in the chart:* $Y = X \circ \phi^{-1}$, *defines a* $d$*-form* $\Omega^k$ *on* $\dot{U}$ *and for every Borel set* $S \subset \dot{U}$,

$$\int_S \Omega^k = \mathrm{E}(M_{u,k}^X(S)).$$

PROOF. Note that a $d$-form is a measure on $\dot{U}$ whose image in each chart is absolutely continuous with respect to Lebesgue measure $\bigwedge_{i=1}^d ds_i$. To prove that (5) defines a $d$-form, it is sufficient to prove that its density with respect to $\bigwedge_{i=1}^d ds_i$ satisfies locally the change-of-variable formula. Let $(U_1, \phi_1)$, $(U_2, \phi_2)$ be two intersecting charts and set

$$U_3 := U_1 \cap U_2; \qquad Y_1 := X \circ \phi_1^{-1}; \qquad Y_2 := X \circ \phi_2^{-1}; \qquad H := \phi_2 \circ \phi_1^{-1}.$$

Denote by $s_i^1$ and $s_i^2$, $i = i, \ldots, d$, the coordinates in each chart. We have

$$\frac{\partial Y_1}{\partial s_i^1} = \sum_{i'} \frac{\partial Y_2}{\partial s_{i'}^2} \frac{\partial H_{i'}}{\partial s_i^1},$$

$$\frac{\partial^2 Y_1}{\partial s_i^1 \, \partial s_j^1} = \sum_{i', j'} \frac{\partial^2 Y_2}{\partial s_{i'}^2 \, \partial s_{j'}^2} \frac{\partial H_{i'}}{\partial s_i^1} \frac{\partial H_{j'}}{\partial s_j^1} + \sum_{i'} \frac{\partial Y_2}{\partial s_{i'}^2} \frac{\partial^2 H_{i'}}{\partial s_i^1 \, \partial s_j^1}.$$

Thus at every point

$$Y_1'(s^1) = (H'(s^1))^T Y_2'(s^2),$$

$$p_{Y_1(s^1), Y_1'(s^1)}(x, 0) = p_{Y_2(s^2), Y_2'(s^2)}(x, 0) |\det(H'(s^1))|^{-1},$$

and at a singular point,

$$Y_1''(s^1) = (H'(s^1))^T Y_2''(s^2) H'(s^1).$$

On the other hand, by the change-of-variable formula,

$$\bigwedge_{i=1}^d ds_i^1 = |\det(H'(s^1))|^{-1} \bigwedge_{i=1}^d ds_i^2.$$

Replacing in the integrand in (5), one checks the desired result.

For the second part again it suffices to prove it locally for an open subset $S$ included in a unique chart. Let $(S, \phi)$ be a chart and let again $Y(s)$ be the process written in this chart. It suffices to check that

$$\begin{aligned}
& \mathrm{E}(M_{u,k}^X(S)) \\
(6) \quad & = \int_{\phi(S)} d\lambda(s) \int_u^{+\infty} dx \, \mathrm{E}(\delta^k(Y''(s))/Y(s) = x, Y'(s) = 0) p_{Y(s), Y'(s)}(x, 0).
\end{aligned}$$

Since $M_{u,k}^X(S)$ is equal to $M_{u,k}^Y\{\phi(S)\}$, we see that the result is a direct consequence of Theorem 2.3.

Even though in the integrand in (5) the product does not depend on the parameterization, each factor does. When the manifold $U$ is equipped with a Riemannian metric it is possible to rewrite (5) as

$$(7) \quad \int_u^{+\infty} dx \; \mathrm{E}(\delta^k(\nabla^2 X(s))/X(s)=x, \nabla X(s)=0) p_{X(s),\nabla X(s)}(x,0) \, \mathrm{Vol},$$

where $\nabla^2 X(s)$ and $\nabla X(s)$ are respectively the Hessian and the gradient read in an orthonormal basis. This formula is close to a formula by Taylor and Adler (2002) for the expected Euler characteristic.

REMARK. One can consider a number of variants of Rice formulae, in which we may be interested in computing the moments of the number of roots of the equation $Z(t)=u$ under some additional conditions. This has been the case in the statement of Theorem 2.3 in which we have given formulae for the first moment of the number of zeroes of $X'$ in which $X$ is bigger than $u$ ($i=2$) and also the real-valued process $X$ has a local maximum ($i=1$).

We just consider below two additional examples of variants that we state here for further reference. We limit the statements to random fields defined on subsets of $\mathcal{R}^d$. Similar statements hold true when the parameter set is a general smooth manifold. Proofs are essentially the same as the previous ones.

VARIANT 1. Assume that $Z_1, Z_2$ are $\mathcal{R}^d$-valued random fields defined on compact subsets $I_1, I_2$ of $\mathcal{R}^d$ and suppose that $(Z_i, I_i)$, $i=1,2$, satisfy the hypotheses of Theorem 2.1 and that for every $s \in I_1$ and $t \in I_2$, the distribution of $(Z_1(s), Z_2(t))$ does not degenerate. Then, for each pair $u_1, u_2 \in \mathcal{R}^d$,

$$\mathrm{E}(N_{u_1}^{Z_1}(I_1) N_{u_2}^{Z_2}(I_2))$$

$$(8) \quad = \int_{I_1 \times I_2} dt_1 \, dt_2 \, \mathrm{E}(|\det(Z_1'(t_1))||\det(Z_2'(t_2))|/Z_1(t_1)=u_1, Z_2(t_2)=u_2)$$

$$\times p_{Z_1(t_1),Z_2(t_2)}(u_1, u_2).$$

VARIANT 2. Let $Z, I$ be as in Theorem 2.1 and let $\xi$ be a real-valued bounded random variable which is measurable with respect to the $\sigma$-algebra generated by the process $Z$. Assume that for each $t \in I$, there exists a continuous Gaussian process $\{Y^t(s) : s \in I\}$, for each $s, t \in I$ a nonrandom function $\alpha^t(s) : \mathcal{R}^d \to \mathcal{R}^d$ and a Borel-measurable function $g : \mathcal{C} \to \mathcal{R}$ where $\mathcal{C}$ is space of real-valued continuous functions on $I$ equipped with the supremum norm, such that:

1. $\xi = g(Y^t(\cdot) + \alpha^t(\cdot)Z(t))$,
2. $Y^t(\cdot)$ and $Z(t)$ are independent,
3. for each $u_0 \in \mathcal{R}$, almost surely the function $u \rightsquigarrow g(Y^t(\cdot) + \alpha^t(\cdot)u)$ is continuous at $u = u_0$.

Then the formula

$$\mathrm{E}(N_u^Z(I)\xi) = \int_I \mathrm{E}(|\det(Z'(t))||\xi/Z(t) = u)p_{Z(t)}(u)\,dt$$

holds true.

We will be particularly interested in the function $\xi = \mathbb{1}_{M_I < v}$ for some $v \in \mathcal{R}$. We will see later on that it satisfies the above conditions under certain hypotheses on the process $Z$.

**3. First derivative, first form.** Our main goals in this and the next section are to prove existence and regularity of the derivative of the function $u \rightsquigarrow F_I(u)$ and, at the same time, that it satisfies some implicit formulae that can be used to provide bounds on it. In the following we assume that $I$ is a $d$-dimensional $\mathcal{C}^\infty$ manifold embedded in $\mathcal{R}^N$, $N \geq d$. $\sigma$ and $\tilde\sigma$ are respectively the geometric measures on $I$ and $\partial I$. Unless explicit statement of the contrary is made, the topology on $I$ will be the relative topology.

In this section we prove formula (10) for $F_I'(u)$—which we call "first form"—which is valid for $\lambda$-almost every $u$, under strong regularity conditions on the paths of the process $X$. In fact, the hypothesis that $X$ is Gaussian is only used in the Rice formula itself and in Lemma 3.1 which gives a bound for the joint density $p_{X(s),X(t),X'(s),X'(t)}$. In both places, one can substitute Gaussianity by appropriate conditions that permit to obtain similar results.

More generally, it is easy to see that inequality (9) is valid under quite general non-Gaussian conditions and implies that we have the following upper bound for the density of the distribution of the random variable $M_I$:

$$\begin{aligned}
F_I'(u) \leq &\int_I \mathrm{E}(\delta^1(X''(t))/X(t) = u, X'(t) = 0)p_{X(t),X'(t)}(u,0)\sigma(dt)\\
&+ \int_{\partial I} \mathrm{E}(\delta^1(\widetilde{X}''(t))/X(t) = u, \widetilde{X}'(t) = 0)p_{X(t),\widetilde{X}'(t)}(u,0)\tilde\sigma(dt),
\end{aligned}$$
(9)

where the function $\delta^1$ has been defined in the statement of Theorem 2.3 and $\widetilde{X}$ denotes the restriction of $X$ to the boundary $\partial I$.

Even for $d = 1$ (one-parameter processes) and $X$ Gaussian and stationary, inequality (9) provides reasonably good upper bounds for $F_I'(u)$ [see Diebolt and Posse (1996) and Azaïs and Wschebor (2001)]. We will see an example for $d = 2$ at the end of this section.

In the next section, we are able to prove that $F_I(u)$ is a $\mathcal{C}^1$ function and that formula (10) can be essentially simplified by getting rid of the conditional expectation, thus obtaining the "second form" for the derivative. This is done under weaker regularity conditions but the assumption that $X$ is Gaussian becomes essential.

DEFINITION 3.1. Let $X : I \to \mathcal{R}$ be a real-valued stochastic process defined on a subset of $\mathcal{R}^d$. We will say that $X$ satisfies condition $(H_k)$, $k$ a positive integer, if the following three conditions hold true:

(i) $X$ is Gaussian;
(ii) a.s. the paths of $X$ are of class $C^k$;
(iii) for any choice of pairwise different values of the parameter $t_1, \ldots, t_n$, the joint distribution of the random variables $X(t_1), \ldots, X(t_n), X'(t_1), \ldots, X'(t_n), \ldots, X^{(k)}(t_1), \ldots, X^{(k)}($
has maximum rank.

The next proposition shows that there exist processes that satisfy $(H_k)$.

PROPOSITION 3.1. Let $X = \{X(t) : t \in \mathcal{R}^d\}$ be a centered stationary Gaussian process having continuous spectral density $f^X$. Assume that $f^X(x) > 0$ for every $x \in \mathcal{R}^d$ and that for any $\alpha > 0$ $f^X(x) \leq C_\alpha \|x\|^{-\alpha}$ holds true for some constant $C_\alpha$ and all $x \in \mathcal{R}^d$. Then, $X$ satisfies $(H_k)$ for every $k = 1, 2, \ldots$.

PROOF. The proof is an adaptation of the proof of a related result for $d = 1$ [Cramér and Leadbetter (1967), page 203]; see Azaïs and Wschebor (2002). □

THEOREM 3.1 (First derivative, first form). Let $X : I \to \mathcal{R}$ be a Gaussian process, $I$ a $\mathcal{C}^\infty$ compact $d$-dimensional manifold. Assume that $X$ verifies $(H_k)$ for every $k = 1, 2, \ldots$.

Then, the function $u \rightsquigarrow F_I(u)$ is absolutely continuous and its Radon–Nikodym derivative is given for almost every $u$ by

$$
F_I'(u) = (-1)^d \int_I \mathrm{E}(\det(X''(t)) \mathbb{1}_{M_I \leq u} / X(t) = u, X'(t) = 0)
$$

$$
\times \, p_{X(t), X'(t)}(u, 0) \sigma(dt)
$$

(10)

$$
+ (-1)^{d-1} \int_{\partial I} \mathrm{E}(\det(\widetilde{X}''(t)) \mathbb{1}_{M_I \leq u} / X(t) = u, \widetilde{X}'(t) = 0)
$$

$$
\times \, p_{X(t), \widetilde{X}'(t)}(u, 0) \widetilde{\sigma}(dt).
$$

PROOF. For $u < v$ and $S$ (resp. $\widetilde{S}$) a subset of $I$ (resp. $\partial I$), let us denote

$$
M_{u,v}(S) = \sharp\{t \in S : u < X(t) \leq v, X'(t) = 0, X''(t) \prec 0\},
$$

$$
\widetilde{M}_{u,v}(\widetilde{S}) = \sharp\{t \in \widetilde{S} : u < X(t) \leq v, \widetilde{X}'(t) = 0, \widetilde{X}''(t) \prec 0\}.
$$

*Step* 1. Let $h > 0$ and consider the increment

$$F_I(u) - F_I(u - h)$$
$$= \mathrm{P}(\{M_I \leq u\} \cap [\{M_{u-h,u}(\dot{I}) \geq 1\} \cup \{\widetilde{M}_{u-h,u}(\partial I) \geq 1\}]).$$

Let us prove that

(11) $$\mathrm{P}(M_{u-h,u}(\dot{I}) \geq 1, \widetilde{M}_{u-h,u}(\partial I) \geq 1) = o(h) \qquad \text{as } h \downarrow 0.$$

In fact, for $\delta > 0$,

(12) $$\mathrm{P}(M_{u-h,u}(\dot{I}) \geq 1, \widetilde{M}_{u-h,u}(\partial I) \geq 1)$$
$$\leq \mathrm{E}(M_{u-h,u}(I_{-\delta}) \widetilde{M}_{u-h,u}(\partial I)) + E(M_{u-h,u}(I \setminus I_{-\delta})).$$

The first term in the right-hand side of (12) can be computed by means of a Rice-type formula, and it can be expressed as

$$\int_{I_{-\delta} \times \partial I} \sigma(dt)\tilde{\sigma}(d\tilde{t}) \iint_{u-h}^{u} dx \, d\tilde{x}$$
$$\times \mathrm{E}(\delta^1(X''(t))\delta^1(\widetilde{X}''(\tilde{t}))/X(t) = x, \widetilde{X}(\tilde{t}) = \tilde{x}, X'(t) = 0, \widetilde{X}'(\tilde{t}) = 0)$$
$$\times p_{X(t),\widetilde{X}(\tilde{t}),X'(t),\widetilde{X}'(\tilde{t})}(x, \tilde{x}, 0, 0),$$

where the function $\delta^1$ has been defined in Theorem 2.3.

Since in this integral $\|t - \tilde{t}\| \geq \delta$, the integrand is bounded and the integral is $O(h^2)$.

For the second term in (12) we apply the Rice formula again. Taking into account that the boundary of $I$ is smooth and compact, we get

$$\mathrm{E}(M_{u-h,u}(I \setminus I_{-\delta}))$$
$$= \int_{I \setminus I_{-\delta}} \sigma(dt) \int_{u-h}^{u} \mathrm{E}(\delta^1(X''(t))/X(t) = x, X'(t) = 0) p_{X(t),X'(t)}(x, 0) \, dx$$
$$\leq (\text{const})h\sigma(I \setminus I_{-\delta}) \leq (\text{const})h\delta,$$

where the constant does not depend on $h$ and $\delta$. Since $\delta > 0$ can be chosen arbitrarily small, (11) follows and we may write as $h \to 0$:

$$F_I(u) - F_I(u - h)$$
$$= \mathrm{P}(M_I \leq u, M_{u-h,u}(\dot{I}) \geq 1) + \mathrm{P}(M_I \leq u, \widetilde{M}_{u-h,u}(\partial I) \geq 1) + o(h).$$

Note that the foregoing argument also implies that $F_I$ is absolutely continuous with respect to Lebesgue measure and that the density is bounded above by the right-hand side of (10). In fact,

$$F_I(u) - F_I(u - h) \leq \mathrm{P}(M_{u-h,u}(\dot{I}) \geq 1) + \mathrm{P}(\widetilde{M}_{u-h,u}(\partial I) \geq 1)$$
$$\leq \mathrm{E}(M_{u-h,u}(\dot{I})) + E(\widetilde{M}_{u-h,u}(\partial I))$$

and it is enough to apply the Rice formula to each one of the expectations on the right-hand side.

The delicate part of the proof consists in showing that we have equality in (10).

*Step* 2. For $g: I \to \mathcal{R}$ we put $\|g\|_\infty = \sup_{t \in I} |g(t)|$ and if $k$ is a nonnegative integer, $\|g\|_{\infty,k} = \sup_{k_1+k_2+\cdots+k_d \le k} \|\partial_{k_1,k_2,\ldots,k_d} g\|_\infty$. For fixed $\gamma > 0$ (to be chosen later on) and $h > 0$, we denote by $E_h = \{\|X\|_{\infty,4} \le h^{-\gamma}\}$. Because of the Landau–Shepp–Fernique inequality [see Landau and Shepp (1970) or Fernique (1975)] there exist positive constants $C_1, C_2$ such that

$$\mathrm{P}(E_h^C) \le C_1 \exp[-C_2 h^{-2\gamma}] = o(h) \qquad \text{as } h \to 0,$$

so that to have (10) it suffices to show that, as $h \to 0$,

$$(13) \qquad \mathrm{E}([M_{u-h,u}(\dot{I}) - \mathbb{1}_{M_{u-h,u}(\dot{I}) \ge 1}] \mathbb{1}_{M_I \le u} \mathbb{1}_{E_h}) = o(h),$$

$$(14) \qquad \mathrm{E}([\widetilde{M}_{u-h,u}(\partial I) - \mathbb{1}_{\widetilde{M}_{u-h,u}(\partial I) \ge 1}] \mathbb{1}_{M_I \le u} \mathbb{1}_{E_h}) = o(h).$$

We prove (13). Equation (14) can be proved in a similar way.

Put $M_{u-h,u} = M_{u-h,u}(\dot{I})$. We have, on applying the Rice formula for the second factorial moment,

$$(15) \quad \begin{aligned} &\mathrm{E}([M_{u-h,u} - \mathbb{1}_{M_{u-h,u} \ge 1}] \mathbb{1}_{M_I \le u} \mathbb{1}_{E_h}) \\ &\le \mathrm{E}(M_{u-h,u}(M_{u-h,u}-1) \mathbb{1}_{E_h}) = \iint_{I \times I} A_{s,t} \sigma(ds) \sigma(dt), \end{aligned}$$

where

$$(16) \quad \begin{aligned} A_{s,t} = \iint_{u-h}^u dx_1 \, dx_2 \\ \times \, \mathrm{E}(|\det(X''(s)) \det(X''(t))| \mathbb{1}_{X''(s) \prec 0, X''(t) \prec 0} \mathbb{1}_{E_h} / X(s) = x_1, \\ X(t) = x_2, X'(s) = 0, X'(t) = 0) \\ \times \, p_{X(s),X(t),X'(s),X'(t)}(x_1, x_2, 0, 0). \end{aligned}$$

Our goal is to prove that $A_{s,t}$ is $o(h)$ as $h \downarrow 0$ uniformly on $s, t$. Note that when $s, t$ vary in a domain of the form $D_\delta := \{t, s \in I : \|t - s\| > \delta\}$ for some $\delta > 0$, then the Gaussian distribution in (16) is nondegenerate and $A_{s,t}$ is bounded by (const)$h^2$, the constant depending on the minimum of the determinant: $\det \mathrm{Var}(X(s), X(t), X'(s), X'(t))$, for $s, t \in D_\delta$.

So it is enough to prove that $A_{s,t} = o(h)$ for $\|t - s\|$ small, and we may assume that $s$ and $t$ are in the same chart $(U, \phi)$. Writing the process in this chart, we may assume that $I$ is a ball or a half ball in $\mathcal{R}^d$. Let $s, t$ be two such points, and define the process $Y = Y^{s,t}$ by $Y(\tau) = X(s + \tau(t-s)); \tau \in [0,1]$. Under the conditioning one has

$$Y(0) = x_1, \qquad Y(1) = x_2, \qquad Y'(0) = Y'(1) = 0,$$

$$Y''(0) = X''(s)[(t-s),(t-s)], \qquad Y''(1) = X''(t)[(t-s),(t-s)].$$

Consider the interpolation polynomial $Q$ of degree 3 such that

$$Q(0) = x_1, \qquad Q(1) = x_2, \qquad Q'(0) = Q'(1) = 0.$$

Check that

$$Q(y) = x_1 + (x_2 - x_1)y^2(3 - 2y), \qquad Q''(0) = -Q''(1) = 6(x_2 - x_1).$$

Denote $Z(\tau) = Y(\tau) - Q(\tau), 0 \le \tau \le 1$. Under the conditioning, one has $Z(0) = Z(1) = Z'(0) = Z'(1) = 0$ and if also the event $E_h$ occurs, an elementary calculation shows that for $0 \le \tau \le 1$,

$$(17) \quad |Z''(\tau)| \le \sup_{\tau \in [0,1]} \frac{|Z^{(4)}(\tau)|}{2!} = \sup_{\tau \in [0,1]} \frac{|Y^{(4)}(\tau)|}{2!} \le (\text{const}) \|t-s\|^4 h^{-\gamma}.$$

On the other hand, check that if $A$ is a positive semidefinite symmetric $d \times d$ real matrix and $v_1$ is a vector of Euclidean norm equal to 1, then the inequality

$$(18) \quad \det(A) \le \langle Av_1, v_1 \rangle \det(B)$$

holds true, where $B$ is the $(d-1) \times (d-1)$ matrix $B = ((\langle Av_j, v_k \rangle))_{j,k=2,\dots,d}$ and $\{v_1, v_2, \dots, v_d\}$ is an orthonormal basis of $R^d$ containing $v_1$.

Assume $X''(s)$ is negative definite, and that the event $E_h$ occurs. We can apply (18) to the matrix $A = -X''(s)$ and the unit vector $v_1 = (t-s)/\|t-s\|$. Note that in that case, the elements of matrix $B$ are of the form $\langle -X''(s)v_j, v_k \rangle$, hence bounded by $(\text{const})h^{-\gamma}$. So,

$$\det[-X''(s)] \le \langle -X''(s)v_1, v_1 \rangle C_d \ h^{-(d-1)\gamma} = C_d[Y''(0)]^{-} \|t-s\|^{-2} h^{-(d-1)\gamma},$$

the constant $C_d$ depending only on the dimension $d$.

Similarly, if $X''(t)$ is negative definite, and the event $E_h$ occurs, then

$$\det[-X''(t)] \le C_d[Y''(1)]^{-} \|t-s\|^{-2} h^{-(d-1)\gamma}.$$

Hence, if $\mathcal{C}$ is the condition $\{X(s) = x_1, X(t) = x_2, X'(s) = 0, X'(t) = 0\}$,

$$\mathrm{E}(|\det(X''(s))\det(X''(t))|\mathbb{1}_{X''(s) \prec 0, X''(t) \prec 0}\mathbb{1}_{E_h}/\mathcal{C})$$

$$\le C_d^2 \ h^{-2(d-1)\gamma} \|t-s\|^{-4} \mathrm{E}([Y''(0)]^{-}[Y''(1)]^{-}\mathbb{1}_{E_h}/\mathcal{C})$$

$$\le C_d^2 \ h^{-2(d-1)\gamma} \|t-s\|^{-4} \mathrm{E}\left(\left[\frac{Y''(0) + Y''(1)}{2}\right]^2 \mathbb{1}_{E_h}/\mathcal{C}\right)$$

$$= C_d^2 \ h^{-2(d-1)\gamma} \|t-s\|^{-4} \mathrm{E}\left(\left[\frac{Z''(0) + Z''(1)}{2}\right]^2 \mathbb{1}_{E_h}/\mathcal{C}\right)$$

$$\le (\text{const}) C_d^2 \ h^{-2d\gamma} \|t-s\|^4.$$

We now turn to the density in (15) using the following lemma which is similar to Lemma 4.3, page 76, in Piterbarg (1996a). The proof is omitted.

LEMMA 3.1.  *For all* $s, t \in I$,

$$(19) \qquad \|t - s\|^{d+3} \, p_{X(s),X(t),X'(s),X'(t)}(0,0,0,0) \leq D,$$

*where* $D$ *is a constant.*

Back to the proof of the theorem, to bound the expression in (15) we use Lemma 3.1 and the bound on the conditional expectation, thus obtaining

$$\mathrm{E}(M_{u-h,u}(M_{u-h,u}-1)\mathbb{1}_{E_h})$$

$$(20) \qquad \leq (\mathrm{const})C_d^2 \, h^{-2d\gamma} D \iint_{I \times I} \|t-s\|^{-d+1} \, ds \, dt \iint_{u-h}^{u} dx_1 \, dx_2$$

$$\leq (\mathrm{const})h^{2-2d\gamma}$$

since the function $(s,t) \rightsquigarrow \|t-s\|^{-d+1}$ is Lebesgue-integrable in $I \times I$. The last constant depends only on the dimension $d$ and the set $I$. Taking $\gamma$ small enough, (13) follows.  □

EXAMPLE.  Let $\{X(s,t)\}$ be a real-valued two-parameter Gaussian, centered stationary isotropic process with covariance $\Gamma$. Assume that $\Gamma(0) = 1$ and that the spectral measure $\mu$ is absolutely continuous with density $\mu(ds, dt) = f(\rho) \, ds \, dt$, $\rho = (s^2 + t^2)^{1/2}$. Assume further that $J_k = \int_0^{+\infty} \rho^k f(\rho) \, d\rho < \infty$, for $1 \leq k \leq 5$. Our aim is to give an explicit upper bound for the density of the probability distribution of $M_I$ where $I$ is the unit disc. Using (9) which is a consequence of Theorem 3.1 and the invariance of the law of the process, we have

$$F_I'(u) \leq \pi \mathrm{E}(\delta^1(X''(0,0))/X(0,0) = u, X'(0,0) = (0,0))$$

$$(21) \qquad \times p_{X(0,0),X'(0,0)}(u,(0,0))$$

$$+ 2\pi \mathrm{E}(\delta^1(\widetilde{X}''(1,0))/X(1,0) = u, \widetilde{X}'(1,0) = 0)p_{X(1,0),\widetilde{X}'(1,0)}(u,0)$$

$$= I_1 + I_2.$$

We denote by $X$, $X'$, $X''$ the value of the different processes at some point $(s,t)$; by $X''_{ss}, X''_{st}, X''_{tt}$ the entries of the matrix $X''$; and by $\varphi$ and $\Phi$ the standard normal density and distribution.

One can easily check that: $X'$ is independent of $X$ and $X''$, and has variance $\pi J_3 I_d$; $X''_{st}$ is independent of $X$, $X'$ $X''_{ss}$ and $X''_{tt}$, and has variance $\frac{\pi}{4} J_5$. Conditionally on $X = u$, the random variables $X''_{ss}$ and $X''_{tt}$ have

$$\text{expectation:} \qquad -\pi J_3;$$

$$\text{variance:} \qquad \frac{3\pi}{4} J_5 - (\pi J_3)^2;$$

$$\text{covariance:} \qquad \frac{\pi}{4} J_5 - (\pi J_3)^2.$$

We obtain

$$I_2 = \sqrt{\frac{2}{J_3}} \varphi(u) \left[ \left( \frac{3\pi}{4} J_5 - (\pi J_3)^2 \right)^{1/2} \varphi(bu) + \pi J_3 u \Phi(bu) \right],$$

with $b = \frac{\pi J_3}{(3\pi/4 J_5 - (\pi J_3)^2)^{1/2}}$. As for $I_1$ we remark that, conditionally on $X = u$, $X''_{ss} + X''_{tt}$ and $X''_{ss} - X''_{tt}$ are independent, so that a direct computation gives

$$
\begin{aligned}
I_1 = \frac{1}{8\pi J_3} \varphi(u) \mathrm{E} \Bigg[ & (\alpha \eta_1 - 2\pi J_3 u)^2 - \frac{\pi J_5}{4} (\eta_2^2 + \eta_3^2) \mathbb{1}_{\{\alpha \eta_1 < 2\pi J_3 u\}} \\
& \times \mathbb{1}_{\{(\alpha \eta_1 - 2\pi J_3 u)^2 - \frac{\pi J_5}{4}(\eta_2^2 + \eta_3^2) > 0\}} \Bigg],
\end{aligned}
$$

(22)

where $\eta_1, \eta_2, \eta_3$ are standard independent normal random variables and $\alpha^2 = 2\pi J_5 - 4\pi^2 J_3^2$. Finally we get

$$
\begin{aligned}
I_1 = \frac{\sqrt{2\pi}}{8\pi J_3} \varphi(u) \int_0^\infty [ & (\alpha^2 + a^2 - c^2 x^2) \Phi(a - cx) \\
& + [2a\alpha - \alpha^2(a - cx)] \varphi(a - cx)] x \varphi(x) \, dx,
\end{aligned}
$$

with $a = 2\pi J_3 u, \ c = \sqrt{\frac{\pi J_5}{4}}$.

**4. First derivative, second form.** We choose, once for this entire section, a finite atlas $\mathcal{A}$ for $I$. Then, to every $t \in I$ it is possible to associate a fixed chart that will be denoted $(U_t, \phi_t)$. When $t \in \partial I$, $\phi_t(U_t)$ can be chosen to be a half ball with $\phi_t(t)$ belonging to the hyperplane limiting this half ball. For $t \in I$, let $V_t$ be an open neighborhood of $t$ whose closure is included in $U_t$ and let $\psi_t$ be a $\mathcal{C}^\infty$ function such that $\psi_t \equiv 1$ on $V_t$; $\psi_t \equiv 0$ on $U_t^c$.

1. For every $t \in \mathring{I}$ and $s \in I$ we define the normalization $n(t, s)$ in the following way:

   (a) For $s \in V_t$, we set "in the chart" $(U_t, \phi_t), n_1(t, s) = \frac{1}{2} \|s - t\|^2$. By "in the chart" we mean that $\|s - t\|$ is in fact $\|\phi_t(t) - \phi_t(s)\|$.
   (b) For general $s$, we set $n(t, s) = \psi_t(s) n_1(t, s) + (1 - \psi_t(s))$.

   Note that in the flat case, when the dimension $d$ of the manifold is equal to the dimension $N$ of the ambient space, the simpler definition $n(t, s) = \frac{1}{2} \|s - t\|^2$ works.

2. For every $t \in \partial I$ and $s \in I$, we set $n_1(t,s) = |(s-t)_N| + \frac{1}{2}\|s-t\|^2$, where $(s-t)_N$ is the normal component of $(s-t)$ with respect to the hyperplane delimiting the half ball $\phi_t(U_t)$. The rest of the definition is the same.

DEFINITION 4.1.    We will say that $f$ is helix-function—or an $h$-function—on $I$ with pole $t \in I$ satisfying hypothesis $(H_{t,k})$, $k$ integer, $k > 1$, if:

(i)  $f$ is a bounded $\mathcal{C}^k$ function on $I \setminus \{t\}$.

(ii)  $\underline{f}(s) := n(t,s)f(s)$ can be prolonged as function of class $\mathcal{C}^k$ on $I$.

DEFINITION 4.2.    In the same way $X$ is called an $h$-process with pole $t \in I$ satisfying hypothesis $(H_{t,k})$, $k$ integer, $k > 1$, if:

(i)  $Z$ is a Gaussian process with $\mathcal{C}^k$ paths on $I \setminus \{t\}$.

(ii)  For $t \in \dot{I}$, $\underline{Z}(s) := n(t,s)Z(s)$ can be prolonged as a process of class $\mathcal{C}^k$ on $I$, with $\underline{Z}(t) = 0$, $\underline{Z}'(t) = 0$, If $s_1,\dots,s_m$ are pairwise different points of $I \setminus \{t\}$, then the distribution of $\underline{Z}^{(2)}(t),\dots,\underline{Z}^{(k)}(t),\underline{Z}(s_1),\dots,\underline{Z}^{(k)}(s_1),\dots,\underline{Z}^{(k)}(s_m)$ does not degenerate.

(iii)  For $t \in \partial I$; $\underline{Z}(s) := n(t,s)Z(s)$ can be prolonged as a process of class $\mathcal{C}^k$ on $I$ with $\underline{Z}(t) = 0$, $\underline{\widetilde{Z}}'(t) = 0$, and if $s_1,\dots,s_m$ are pairwise different points of $I \setminus \{t\}$, then the distribution of $\underline{Z}'_N(t),\underline{Z}^{(2)}(t),\dots,\underline{Z}^{(k)}(t),\underline{Z}(s_1),\dots,\underline{Z}^{(k)}(s_1),\dots,\underline{Z}^{(k)}(s_m)$ does not degenerate. $\underline{Z}'_N(t)$ is the derivative normal to the boundary of $I$ at $t$.

We use the terms "$h$-function" and "$h$-process" since the function and the paths of the process need not extend to a continuous function at the point $t$. However, the definition implies the existence of radial limits at $t$. So the process may take the form of a helix around $t$.

LEMMA 4.1.    Let $X$ be a process satisfying $(H_k, k \geq 2)$, and let $f$ be a $\mathcal{C}^k$ function $I \to \mathcal{R}$.

(a) For $t \in \dot{I}$, set for $s \in I, s \neq t$,

$$X(s) = a_s^t X(t) + \langle b_s^t, X'(t)\rangle + n(t,s)X^t(s),$$

where $a_s^t$ and $b_s^t$ are the regression coefficients.

In the same way, set

$$f(s) = a_s^t f(t) + \langle b_s^t, f'(t)\rangle + n(t,s)f^t(s),$$

using the regression coefficients associated to $X$.

(b) For $t \in \partial I$, $s \in T, s \neq t$, set

$$X(s) = \tilde{a}_s^t X(t) + \langle \tilde{b}_s^t, \widetilde{X}'(t)\rangle + n(t,s)X^t(s)$$

and

$$f(s) = \tilde{a}_s^t f(t) + \langle \tilde{b}_s^t, \tilde{f}'(t)\rangle + n(t,s)f^t(s).$$

*Then $s \rightsquigarrow X^t(s)$ and $s \rightsquigarrow f^t(s)$ are, respectively, an h-process and an h-function with pole $t$ satisfying $H_{t,k}$.*

PROOF. We give the proof in the case $t \in \dot{I}$, the other one being similar. In fact, the quantity denoted by $\underline{X}^t(s)$ is just $X(s) - a_s^t X(t) - \langle b_s^t, X'(t) \rangle$. On $L^2(\Omega, P)$, let $\Pi$ be the projector on the orthogonal complement to the subspace generated by $X(t), X'(t)$. Using a Taylor expansion,

$$X(s) = X(t) + \langle (s - t), X'(t) \rangle$$
$$+ \|t - s\|^2 \int_0^1 X''((1 - \alpha)t + \alpha s)[v, v](1 - \alpha) \, d\alpha,$$

with $v = \frac{s-t}{\|s-t\|}$. This implies that

$$(23) \qquad \underline{X}^t(s) = \Pi \bigg[ \|t - s\|^2 \int_0^1 X''((1 - \alpha)t + \alpha s)[v, v](1 - \alpha) \, d\alpha \bigg],$$

which gives the result due to the nondegeneracy condition. $\quad\square$

We state now an extension of Ylvisaker's (1968) theorem on the regularity of the distribution of the maximum of a Gaussian process which we will use in the proof of Theorem 4.2 and which might have some interest in itself.

THEOREM 4.1. *Let $Z : T \to \mathcal{R}$ be a Gaussian separable process on some parameter set $T$ and denote by $M^Z = \sup_{t \in T} Z(t)$ which is (a random variable) taking values in $\mathcal{R} \cup \{+\infty\}$. Assume that there exists $\sigma_0 > 0$, $m_- > -\infty$ such that*

$$m(t) = \mathrm{E}(Z_t) \geq m_-, \qquad \sigma^2(t) = \mathrm{Var}(Z_t) \geq \sigma_0^2$$

*for every $t \in T$. Then the distribution of the random variable $M^Z$ is the sum of an atom at $+\infty$ and a—possibly defective—probability measure on $\mathcal{R}$ which has a locally bounded density.*

PROOF. Suppose first that $X : T \to \mathcal{R}$ is a Gaussian separable process satisfying $\mathrm{Var}(X_t) = 1$, $\mathrm{E}(X_t) \geq 0$, for every $t \in T$. A close look at Ylvisaker's (1968) proof shows that the distribution of the supremum $M^X$ has a density $p_{M^X}$ that satisfies

$$(24) \qquad p_{M^X}(u) \leq \psi(u) = \frac{\exp(-u^2/2)}{\int_u^\infty \exp(-v^2/2) \, dv} \qquad \text{for every } u \in \mathcal{R}.$$

Let now $Z$ satisfy the hypotheses of the theorem. For given $a, b \in \mathcal{R}$, $a < b$, choose $A \in \mathcal{R}^+$ so that $|a| < A$ and consider the process

$$X(t) = \frac{Z(t) - a}{\sigma(t)} + \frac{|m_-| + A}{\sigma_0}.$$

Clearly, for every $t \in T$,

$$\mathrm{E}(X(t)) = \frac{m(t) - a}{\sigma(t)} + \frac{|m_-| + A}{\sigma_0} \geq -\frac{|m_-| + |a|}{\sigma_0} + \frac{|m_-| + A}{\sigma_0} \geq 0,$$

and $\mathrm{Var}(X(t)) = 1$. So that (24) holds for the process $X$.

On the other hand, the statement follows from the inclusion:

$$\{a < M^Z \leq b\} \subset \left\{ \frac{|m_-| + A}{\sigma_0} < M^X \leq \frac{|m_-| + A}{\sigma_0} + \frac{b - a}{\sigma_0} \right\},$$

which implies

$$\mathrm{P}\{a < M^Z \leq b\} \leq \int_{(|m_-|+A)/\sigma_0}^{(|m_-|+A)/\sigma_0 + (b-a)/\sigma_0} \psi(u)\, du$$

$$= \int_a^b \frac{1}{\sigma_0} \psi\left( \frac{v - a + |m_-| + A}{\sigma_0} \right) dv. \qquad \square$$

Set now $\beta(t) \equiv 1$. The key point is that, due to regression formulae, under the condition $\{X(t) = u, X'(t) = 0\}$ the event

$$A_u(X, \beta) := \{X(s) \leq u, \forall s \in I\}$$

coincides with the event

$$A_u(X^t, \beta^t) := \{X^t(s) \leq \beta^t(s) u, \forall s \in I \setminus \{t\}\},$$

where $X^t$ and $\beta^t$ are the h-process and the h-function defined in Lemma 4.1.

THEOREM 4.2 (First derivative, second form). *Let $X : I \to \mathcal{R}$ be a Gaussian process, $I$ a $C^\infty$ compact manifold contained in $\mathcal{R}^d$. Assume that $X$ has paths of class $C^2$ and for $s \neq t$ the triplet $(X(s), X(t), X'(t))$ in $\mathcal{R} \times \mathcal{R} \times \mathcal{R}^d$ has a nondegenerate distribution. Then, the result of Theorem 3.1 is valid, the derivative $F_I'(u)$ given by relation (10) can be written as*

$$F_I'(u) = (-1)^d \int_I \mathrm{E}[\det(\underline{X}'^t(t) - \underline{\beta}'^t(t)u) \mathbb{1}_{A_u(X^t, \beta^t)}]$$

$$\times p_{X(t), X'(t)}(u, 0) \sigma(dt)$$

$$(25) \qquad + (-1)^{d-1} \int_{\partial I} \mathrm{E}[\det(\underline{\widetilde{X}}'^t(t) - \underline{\widetilde{\beta}}'^t(t))u \mathbb{1}_{A_u(X^t, \beta^t)}]$$

$$\times p_{X(t), \widetilde{X}'(t)}(u, 0) \tilde{\sigma}(dt)$$

*and this expression is continuous as a function of $u$.*

The notation $\widetilde{\underline{X}}^{t''}(t)$ should be understood in the sense that we first define $\underline{X}^t$ and then calculate its second derivative along $\partial I$.

PROOF OF THEOREM 4.2. As a first step, assume that the process $X$ satisfies the hypotheses of Theorem 3.1, which are stronger that those in the present theorem.

We prove that the first term in (10) can be rewritten as the first term in (25). One can proceed in a similar way with the second term, mutatis mutandis. For that purpose, use the remark just before the statement of Theorem 4.2 and the fact that under the condition $\{X(t) = u, X'(t) = 0\}$, $X''(t)$ is equal to $\underline{X}^{t''}(t) - \underline{\beta}^{t''}(t)u$.

Replacing in the conditional expectation in (10) and on account of the Gaussianity of the process, we get rid of the conditioning and obtain the first term in (25). We now study the continuity of $u \rightsquigarrow F_I'(u)$. The variable $u$ appears at three locations:

(i) in the density $p_{X(t),X'(t)}(u, 0)$, which is clearly continuous,
(ii) in

$$\mathrm{E}[\det(\underline{X}^{t''}(t) - \underline{\beta}^{t''}(t)u)\mathbb{1}_{A_u(X^t, \beta^t)}],$$

where it occurs twice: in the first factor and in the indicator function.

Due to the integrability of the supremum of bounded Gaussian processes, it is easy to prove that this expression is continuous as a function of the first $u$.

As for the $u$ in the indicator function, set

$$(26) \qquad \xi_v := \det(\underline{X}^{t''}(t) - \underline{\beta}^{t''}(t)v)$$

and, for $h > 0$, consider the quantity $\mathrm{E}[\xi_v \mathbb{1}_{A_u(X^t, \beta^t)}] - \mathrm{E}[\xi_v \mathbb{1}_{A_{u-h}(X^t, \beta^t)}]$, which is equal to

$$(27) \qquad \mathrm{E}[\xi_v \mathbb{1}_{A_u(X^t, \beta^t) \setminus A_{u-h}(X^t, \beta^t)}] - \mathrm{E}[\xi_v \mathbb{1}_{A_{u-h}(X^t, \beta^t) \setminus A_u(X^t, \beta^t)}].$$

Apply Schwarz's inequality to the first term in (27):

$$\mathrm{E}[\xi_v \mathbb{1}_{A_u(X^t, \beta^t) \setminus A_{u-h}(X^t, \beta^t)}] \leq [\mathrm{E}(\xi_v^2)\mathrm{P}\{A_u(X^t, \beta^t) \setminus A_{u-h}(X^t, \beta^t)\}]^{1/2}.$$

The event $A_u(X^t, \beta^t) \setminus A_{u-h}(X^t, \beta^t)$ can be described as

$$\forall s \in I \setminus \{t\} : X^t(s) - \beta^t(s)u \leq 0; \qquad \exists s_0 \in I \setminus \{t\} : X^t(s_0) - \beta^t(s_0)(u-h) > 0.$$

This implies that $\beta^t(s_0) > 0$ and that $-\|\beta^t\|_\infty h \leq \sup_{s \in I \setminus \{t\}} X^t(s) - \beta^t(s)u \leq 0$. Now, observe that our improved version of Ylvisaker's theorem (Theorem 4.1) applies to the process $s \rightsquigarrow X^t(s) - \beta^t(s)u$ defined on $I \setminus \{t\}$. This implies that the first term in (27) tends to zero as $h \downarrow 0$. An analogous argument

applies to the second term. Finally, the continuity of $F_I'(u)$ follows from the fact that one can pass to the limit under the integral sign in (25).

To complete the proof we still have to show that the added hypotheses are in fact unnecessary for the validity of the conclusion. Suppose now that the process $X$ satisfies only the hypotheses of the theorem and define

$$(28) \qquad X^\varepsilon(t) = Z_\varepsilon(t) + \varepsilon Y(t),$$

where for each $\varepsilon > 0$, $Z_\varepsilon$ is a real-valued Gaussian process defined on $I$, measurable with respect to the $\sigma$-algebra generated by $\{X(t): t \in I\}$, possessing $\mathcal{C}^\infty$ paths and such that almost surely $Z_\varepsilon(t)$, $Z_\varepsilon'(t)$, $Z_\varepsilon''(t)$ converge uniformly on $I$ to $X(t), X'(t), X''(t)$, respectively, as $\varepsilon \downarrow 0$. One standard form to construct such an approximation process $Z_\varepsilon$ is to use a $\mathcal{C}^\infty$ partition of the unity on $I$ and to approximate locally the composition of a chart with the function $X$ by means of a convolution with a $\mathcal{C}^\infty$ kernel.

In (28), $Y$ denotes the restriction to $I$ of a Gaussian centered stationary process satisfying the hypotheses of Proposition 3.1, defined on $\mathcal{R}^N$, and independent of $X$. Clearly $X^\varepsilon$ satisfies condition $(H_k)$ for every $k$, since it has $\mathcal{C}^\infty$ paths and the independence of both terms in (28) ensures that $X^\varepsilon$ inherits from $Y$ the nondegeneracy condition in Definition 3.1. So, if $M_I^\varepsilon = \max_{t \in I} X^\varepsilon(t)$ and $F_I^\varepsilon(u) = \mathrm{P}\{M_I^\varepsilon \leq u\}$, one has

$$
\begin{aligned}
(29) \qquad F_I^{\varepsilon\prime}(u) = {}& (-1)^d \int_I \mathrm{E}[\det(\underline{X}^{\varepsilon t\prime\prime}(t) - \underline{\beta}^{\varepsilon t\prime\prime}(t)u)\mathbb{1}_{A_u(X^{\varepsilon t}, \beta^{\varepsilon, t})}] \\
& \qquad \times p_{X^\varepsilon(t), X^{\varepsilon\prime}(t)}(u, 0)\sigma(dt) \\
& + (-1)^{d-1} \int_{\partial I} \mathrm{E}[\det(\underline{\widetilde{X}}^{\varepsilon t\prime\prime}(t) - \underline{\widetilde{\beta}}^{\varepsilon t\prime\prime}(t)u)\mathbb{1}_{A_u(X^{\varepsilon t}, \beta^{\varepsilon t})}] \\
& \qquad \times p_{X^\varepsilon(t), \widetilde{X}^{\varepsilon\prime}(t)}(u, 0)\tilde{\sigma}(dt).
\end{aligned}
$$

We want to pass to the limit as $\varepsilon \downarrow 0$ in (29). We prove that the right-hand side is bounded if $\varepsilon$ is small enough and converges to a continuous function of $u$ as $\varepsilon \downarrow 0$. Since $M_I^\varepsilon \to M_I$, this implies that the limit is continuous and coincides with $F_I'(u)$ by a standard argument on convergence of densities. We consider only the first term in (29); the second is similar.

The convergence of $X_\varepsilon$ and its first and second derivative, together with the nondegeneracy hypothesis, imply that uniformly on $t \in I$, as $\varepsilon \downarrow 0$, $p_{X^\varepsilon(t), X^{\varepsilon\prime}(t)}(u, 0) \to p_{X(t), X'(t)}(u, 0)$. The same kind of argument can be used for $\det(\underline{X}^{\varepsilon t\prime\prime}(t) - \underline{\beta}^{\varepsilon t\prime\prime}(t)u)$, on account of the form of the regression coefficients and the definitions of $\underline{X}^t$ and $\underline{\beta}^t$. The only difficulty is to prove that, for fixed $u$,

$$(30) \qquad \mathrm{P}\{C_\varepsilon \Delta C\} \to 0 \qquad \text{as } \varepsilon \downarrow 0,$$

where $C_\varepsilon = A_u(X^{\varepsilon t}, \beta^{\varepsilon t})$, $C = A_u(X^t, \beta^t)$.

We prove that

$$\text{(31)} \qquad \text{a.s. } \mathbb{1}_{C_\varepsilon} \to \mathbb{1}_C \qquad \text{as } \varepsilon \downarrow 0,$$

which implies (30). First of all, note that the event

$$L = \left\{ \sup_{s \in I \setminus \{t\}} (X^t(s) - \beta^t(s)u) = 0 \right\}$$

has zero probability, as already mentioned. Second, from the definition of $X^t(s)$ and the hypothesis, it follows that, as $\varepsilon \downarrow 0$, $X^{\varepsilon, t}(s), \beta^{\varepsilon, t}(s)$ converge to $X^t(s), \beta^t(s)$ uniformly on $I \setminus \{t\}$. Now, if $\omega \notin C$, there exists $\bar{s} = \bar{s}(\omega) \in I \setminus \{t\}$ such that $X^t(\bar{s}) - \beta^t(\bar{s})u > 0$ and for $\varepsilon > 0$ small enough, one has $X^{\varepsilon t}(\bar{s}) - \beta^{\varepsilon t}(\bar{s})u > 0$, which implies that $\omega \notin C_\varepsilon$.

On the other hand, let $\omega \in C \setminus L$. This implies that

$$\sup_{s \in I \setminus \{t\}} (X^t(s) - \beta^t(s)u) < 0.$$

From the above-mentioned uniform convergence, it follows that if $\varepsilon > 0$ is small enough, then $\sup_{s \in I \setminus \{t\}} (X^{\varepsilon t}(s) - \beta^{\varepsilon t}(s)u) < 0$, hence $\omega \in C_\varepsilon$. Equation (31) follows.

So, we have proved that the limit as $\varepsilon \downarrow 0$ of the first term in (29) is equal to the first term in (25).

It remains only to prove that the first term in (25) is a continuous function of $u$. For this purpose, it suffices to show that the function $u \rightsquigarrow \mathrm{P}\{A_u(X^t, \beta^t)\}$ is continuous. This is a consequence of the inequality

$$|\mathrm{P}\{A_{u+h}(X^t, \beta^t)\} - \mathrm{P}\{A_u(X^t, \beta^t)\}|$$
$$\leq \mathrm{P}\left\{ \left| \sup_{s \in I \setminus \{t\}} (X^t(s) - \beta^t(s)u) \right| \leq |h| \sup_{s \in I \setminus \{t\}} |\beta^t(s)| \right\}$$

and of Theorem 4.1, applied once again to the process $s \rightsquigarrow X^t(s) - \beta^t(s)u$ defined on $I \setminus \{t\}$. $\quad\square$

## 5. Asymptotic expansion of $F'(u)$ for large $u$.

COROLLARY 5.1. *Suppose that the process $X$ satisfies the conditions of Theorem 4.2 and that in addition $\mathrm{E}(X_t) = 0$ and $\mathrm{Var}(X_t) = 1$.*

*Then, as $u \to +\infty$, $F'(u)$ is equivalent to*

$$\text{(32)} \qquad \frac{u^d}{(2\pi)^{(d+1)/2}} e^{-u^2/2} \int_I (\det(\Lambda(t)))^{1/2} \, dt,$$

*where $\Lambda(t)$ is the variance–covariance matrix of $X'(t)$.*

Note that (32) is in fact the derivative of the bound for the distribution function that can be obtained by Rice's method [Azaïs and Delmas (2002)] or by the expected Euler characteristic method [Taylor, Takemura and Adler (2004)].

PROOF OF COROLLARY 5.1.    Set $r(s,t) := \mathrm{E}(X(s), X(t))$, and for $i, j = 1, d,$

$$r_{i;}(s,t) := \frac{\partial}{\partial s_i} r(s,t),$$

$$r_{ij;}(s,t) := \frac{\partial^2}{\partial s_i \partial s_j} r(s,t), \qquad r_{i;j}(s,t) := \frac{\partial^2}{\partial s_i \partial t_j} r(s,t).$$

For every $t, i$ and $j$, $r_{i;}(t,t) = 0$, $\Lambda_{ij}(t) = r_{i;j}(t,t) = -r_{ij;}(t,t)$. Thus $X(t)$ and $X'(t)$ are independent. Regression formulae imply that $a_s^t = r(s,t)$, $\beta^t(s) = \frac{1-r(t,s)}{n(s,t)}$. This implies that $\underline{\beta}^t(t) = \Lambda(t)$ and that the possible limits values of $\beta^t(s)$ as $s \to t$ are in the set $\{v^T \Lambda(t) v : v \in S^{d-1}\}$. Due to the nondegeneracy condition these quantities are minorized by a positive constant. On the other hand, for $s \neq t$, $\beta^t(s) > 0$. This shows that for every $t \in I$ one has $\inf_{s \in I} \beta^t(s) > 0$. Since for every $t \in I$ the process $X^t$ is bounded, it follows that a.s. $\mathbb{1}_{A_u(X^t, \beta^t)} \to 1$ as $u \to +\infty$. Also

$$\det(\underline{X}^{t\prime\prime}(t) - \underline{\beta}^{t\prime\prime}(t)u) \simeq (-1)^d \det(\Lambda(t))u^d.$$

Dominated convergence shows that the first term in (25) is equivalent to

$$\int_I u^d \det(\Lambda^t)(2\pi)^{-1/2} e^{-u^2/2} (2\pi)^{-d/2} (\det(\Lambda^t))^{-1/2} \, dt$$

$$= \frac{u^d}{(2\pi)^{(d+1)/2}} e^{-u^2/2} \int_I (\det(\Lambda^t))^{1/2} \, dt.$$

The same kind of argument shows that the second term is $O(u^{d-1} e^{-u^2/2})$, which completes the proof.  $\square$

**Acknowledgment.**   We thank an anonymous referee for very carefully reading the first version of this work and for very valuable suggestions.

LABORATOIRE DE STATISTIQUE
ET PROBABILITÉS
UMR-CNRS C5583
UNIVERSITÉ PAUL SABATIER
118 ROUTE DE NARBONNE
31062 TOULOUSE CEDEX 4
FRANCE
E-MAIL: azais@cict.fr

CENTRO DE MATEMÁTICA
FACULTAD DE CIENCIAS
UNIVERSIDAD DE LA REPÚBLICA
CALLE IGUA 4225
11400 MONTEVIDEO
URUGUAY
E-MAIL: wscheb@fcien.edu.uy